\documentclass[review,11pt,authoryear]{elsarticle}
\usepackage[singlelinecheck=false,labelsep=period,font=footnotesize]{caption}
\usepackage{comment}
\usepackage{dsfont}
\usepackage{float}
\usepackage{multicol}
\usepackage{multirow}
\usepackage{hyperref}
\usepackage{graphicx}
\usepackage{amsfonts}
\usepackage{amsmath}
\usepackage{color}
\usepackage{epstopdf}
\usepackage{subcaption}
\usepackage{eurosym}
\usepackage{enumitem}
\usepackage{bbm}
\usepackage{tikz,pgfplots,pgfplotstable,booktabs}
\usetikzlibrary{calc}
\pgfplotsset{compat=1.8}
\usepgfplotslibrary{statistics}
\usepackage{balance}

\DeclareMathOperator*{\argmin}{arg\,min} 
\usepackage{lineno}
\usepackage{mathtools}
\DeclarePairedDelimiter{\ceil}{\lceil}{\rceil}
\allowdisplaybreaks

\tikzset{sin v source/.style={
  circle,
  draw,
  append after command={
    \pgfextra{
    \draw
      ($(\tikzlastnode.center)!0.5!(\tikzlastnode.west)$)
       arc[start angle=180,end angle=0,radius=0.425ex] 
      (\tikzlastnode.center)
       arc[start angle=180,end angle=360,radius=0.425ex]
      ($(\tikzlastnode.center)!0.5!(\tikzlastnode.east)$) 
    ;
    }
  },
  scale=1.5,
 }
}

\begin{document}

\begin{frontmatter}
    \title{Prescribing net demand for two-stage electricity generation scheduling}
    \author{J. M. Morales\corref{mycorrespondingauthor}}  
    \cortext[mycorrespondingauthor]{Corresponding author}
    \ead{juan.morales@uma.es}
    \author{M. A. Mu\~noz}
    \ead{miguelangeljmd@uma.es}
    \author{S. Pineda}
    \ead{spinedamorente@gmail.com}
    \address{OASYS Group, Ada Byron Research Building,
Arquitecto Francisco Pe\~nalosa St., 18, 29010,
University of Malaga, Malaga, Spain}
    
    \begin{abstract}
    We consider a two-stage generation scheduling problem comprising a forward dispatch and a real-time re-dispatch. The former must be conducted facing an uncertain net demand that includes non-dispatchable electricity consumption and renewable power generation. The latter copes with the plausible deviations with respect to the forward schedule by making use of balancing power during the actual operation of the system. Standard industry practice deals with the uncertain net demand in the forward stage by replacing it with a good estimate of its conditional expectation (usually referred to as a \emph{point forecast}), so as to minimize the need for balancing power in real time. However, it is well known that the cost structure of a power system is highly asymmetric and dependent on its operating point, with the result that minimizing the amount of power imbalances is not necessarily aligned with minimizing operating costs. In this paper, we propose a bilevel program to construct, from the available historical data, a \emph{prescription} of the {\color{black} net demand} that does account for the power system's cost asymmetry. Furthermore, to accommodate the strong dependence of this cost on the power system's operating point, we use clustering to tailor the proposed prescription to the foreseen net-demand regime. By way of an illustrative example and a more realistic case study based on the European power system, we show that our approach leads to substantial cost savings compared to the customary way of doing.
    \end{abstract}
    
    \begin{keyword}
    Smart predict \sep Net demand prescription \sep Two-stage power generation scheduling \sep
    Data-driven optimization  
    \end{keyword}
\end{frontmatter}

\section{Introduction}

Many decision-making processes under uncertainty can be modeled by optimization problems where some of the input parameters are not perfectly known. The field of Optimization under Uncertainty focuses on developing tools to tackle these problems depending on the knowledge of those parameters that the decision maker actually has. For example, if these parameters can be modeled reasonably well as random variables following certain probability distributions, then the decision maker should probably resort to stochastic programming techniques~\citep{birge2011introduction}. In contrast, if all the decision maker knows about said parameters is their range of variation or support, then she should rather opt for robust optimization methods instead~\citep{ben2009robust}.

In the realm of power system operations, there is a vast literature on \textcolor{black}{operations research} methods, models, and algorithms for power dispatch that rely on stochastic programming or robust optimization or hybrids of both. The richness of this literature makes it materially impossible and pointless to embrace it all in this paper. Instead, we refer the reader to monograph~\cite{morales2013integrating} and references therein for examples of power systems operating problems based on stochastic programming, to the seminal work~\cite{bertsimas2012adaptive} on the application of robust optimization for unit commitment and generation scheduling, and to the recent contribution~\cite{dvorkin2019chance} on a distributionally robust chance-constrained electricity market.

Despite the firm and promising advances in Optimization under Uncertainty, still one of the most widely extended practices in decision making is to replace the unknown parameter with a sensible value or estimate, some sort of ``the most likely value'' that the parameter can take on. A natural candidate to play that role is the expected value of the parameter. Thus,  the decision maker can, in addition, exploit all the powerful tools that the disciplines of Statistics, Forecasting and Machine Learning have developed for decades to estimate that expected value \emph{conditional} on all the information the decision maker has available at the moment the decision must be made. The adherence of the power sector to this strategy is particularly notorious, essentially because it is argued to be simpler, more transparent, computationally cheaper and easily accepted by the different stakeholders (see, e.g., \cite{morales2014electricity, wang2015real, morales2017inefficiency, kazempour2018stochastic} for further details on this issue). However, \textcolor{black}{operations researchers} have repeatedly shown that this strategy results in suboptimal decisions in general, because the conditional expected value of the parameter ignores the impact of the parameter's uncertainty on the decision's value~\citep{birge2011introduction}.

Against this background, research efforts have been recently placed on finding a compromise solution. Intuitively, the idea is still to replace the uncertain parameter but with a point estimate (generally different from the parameter's conditional expectation) that is purposely computed to result in the optimal or a nearly optimal decision in view of the parameter's uncertainty. This alternative point estimate is usually called a \emph{prescription} and the term \emph{smart predict} has been recently coined to refer to this midway solution strategy. In this line, we find a number of research works, e.g., \cite{Donti2017, munoz2021bilevel, Elmachtoub2021, balghiti2019generalization}. In particular,  the authors in~\cite{Donti2017} propose a heuristic gradient-based procedure to produce estimates of uncertain parameters in optimization problems based on the objective of the task for which these estimates will be used. In \cite{munoz2021bilevel}, instead, they introduce a bilevel programming framework to the same end. In \cite{Elmachtoub2021}, they deal with linear programs with an uncertain cost vector, for which they develop a tailored convex loss function to compute an estimate of the cost vector that accounts for the underlying linear optimization problem. Finally, the work in~\cite{balghiti2019generalization} is a continuation of that in~\cite{Elmachtoub2021} {\color{black}(first released as an arXiv preprint in 2017)}, where the authors provide bounds on how well a certain method to predict the cost vector from training data generalizes out of sample. 

In the field of power systems, it has also been shown that, by smartly tuning the input parameters of current operational and procedures, these can mimic the performance of their stochastic-programming-based counterparts to a large extent. For instance, in~\cite{morales2014electricity}, they propose a bilevel programming model to compute the amount of (uncertain) renewable power generation that must be considered in a forward electricity market to maximize the short-run market efficiency. In the same vein, the authors in~\cite{dvorkin2018setting} show that, {\color{black} by means of bilevel programming too, the reserve requirements in an European-style two-stage electricity market can be set so that the market can be almost as cost-efficient as the ideal two-stage electricity market run by a full stochastic programming approach}. {\color{black} Our proposal is in the same spirit, but we make use of bilevel programming not to directly prescribe renewable production or reserve capacity, but to determine straightforward \emph{rules} to improve the value of the (net-demand) forecasts currently employed by system operators. In any case, the works of ~\cite{morales2014electricity} and ~\cite{dvorkin2018setting}} reveal that the power sector can highly benefit from the aforementioned smart-predict strategy. Actually, in~\cite{Donti2017}, they apply it to power generation and grid-scale electricity battery operation, and in~\cite{munoz2021bilevel} to the offering problem of a thermal power producer competing strategically in an electricity market. In \cite{carriere2019integrated} and \cite{munoz2020feature}, they focus instead on renewable energy producers, for which they propose different smart-predict strategies for energy trading. Lastly, the authors in~\cite{garcia2021} use a bilevel programming framework similar to that proposed in~\cite{munoz2021bilevel} whereby they train several autoregressive models to estimate the uncertain demand and the size of the energy reserves in a joint reserve allocation and energy dispatch problem. 

Within this context, and given a two-stage power scheduling problem, our contributions are the following:
\begin{itemize}
    \item We propose a bilevel program to learn the value of the net demand (i.e., that obtained by subtracting the weather-driven renewable power production from the non-dispatchable power load) that the scheduling problem must consider so that the power system operating costs are minimized in expectation. This value  (that is, the prescription) is, in general, different from the conditional mean of the net demand that is currently employed in standard industry practice. 
    \item We show that our method can be easily used to upgrade \emph{any} given point forecast of the net demand into a prescription that accounts for the power system's cost asymmetry.
    \item {\color{black} Motivated by the fact that a power system usually features distinct net-demand regimes,} we introduce a data clustering and partitioning strategy that, on the one hand, increases the performance of the  prescription of the net demand that our approach produces (by making it dependent on the foreseen net-demand regime), and, on the other, substantially decreases the computational effort to solve the associated bilevel estimation problem.
    \item We evaluate the economic benefits that our approach achieves through an out-of-sample test on a stylized version of the European power system that makes use of real data, in particular, of actual and day-ahead predicted net-demand values downloaded from the ENTSO-e Transaparency Platform~\citep{Entsoe}.
\end{itemize}

The rest of this paper is organized as follows. Section~\ref{sec:PD} describes the two-stage power generation scheduling problem we consider throughout our paper, motivates the ultimate goal of our work, and formulates the generic bilevel program we use to construct, from the available historical data, prescriptions of the net demand intended to minimize the expected system operation costs. In Section~\ref{sec:CD}, we resort to various simplifying assumptions to guarantee that the globally optimal solution of the bilevel program can be obtained by solving a mixed-integer linear program (MILP). The potential of the net-demand prescriptions provided by this MILP is then illustrated, discussed and justified using a small power system in Section~\ref{sec:Ex}. In Section~\ref{sec:Part}, we introduce a procedure based on data clustering and partitioning to enhance the value of our prescription and to speed up the solution of the mixed-integer prescription problem. A case study based on the European power system is used in Section~\ref{sec:CaseStudy} to investigate the benefits of our approach on real data. Finally, Section~\ref{sec:Conclu} concludes the paper with some final remarks.

\section{Problem description}\label{sec:PD}
%

We consider a two-stage generation scheduling problem consisting of a forward power dispatch and a real-time balancing phase. The forward dispatch is decided some time prior to the actual delivery of energy, for instance, from 15 minutes to 36 hours in advance. \textcolor{black}{The forward dispatch is required to decide the output levels of inflexible generating units (e.g. nuclear-based power plants) that need advance planning due to technical limits such as ramping constraints or minimum times.} The real-time balancing phase processes the energy imbalances with respect to the forward production schedule. \textcolor{black}{The real-time balancing rescheduling aims at adjusting the output of flexible generating units (such as gas-based power plants) that can quickly deviate from the forward schedule to ensure the equilibrium between electricity production and consumption. More details about the two-stage generation scheduling problem used in this paper can be found in \cite{morales2013integrating}}

In our framework, the forward stage first determines the power dispatch that minimizes the anticipated electricity production costs by solving
\begin{subequations}
\begin{align}
    \min_{\mathbf{z}^{\rm{F}}} & \enskip f^{\text{obj}}(\mathbf{z}^{\rm{F}}, \hat{\mathbf{y}}) \label{eq:generic_forward_obj} \\
    \text{s.t.} & \enskip f^{\text{gen}}(\mathbf{z}^{\rm{F}}, \hat{\mathbf{y}}) \leq 0 \label{eq:generic_forward_gen} \\
    \phantom{s.t.} & \enskip f^{\text{net}}(\mathbf{z}^{\rm{F}}, \hat{\mathbf{y}}) \leq 0, \label{eq:generic_forward_net}
\end{align}\label{eq:generic_forward}
\end{subequations}
\noindent where $\mathbf{z}^{\rm{F}}$ denotes the scheduling decisions (such as the dispatch of thermal generating units and the consumption of dispatchable loads) and $\hat{\mathbf{y}}$ represents a point or single-value estimate of the uncertain parameters (in our case, the net demand, given as the difference between the non-dispatchable load and the weather-driven renewable power generation). Objective function~\eqref{eq:generic_forward_obj} determines the forward system operating costs, while the generic constraints~\eqref{eq:generic_forward_gen} and \eqref{eq:generic_forward_net} enforce the technical limitations of the generating units and the network, respectively. 

Since the net demand is uncertain, the scheduling decisions that results from the forward problem~\eqref{eq:generic_forward} are to be adjusted during the real-time operation of the power system to accommodate the \emph{eventual} net-demand value that is realized. This adjustment is conducted through the following generic real-time balancing problem:
\begin{subequations}
\begin{align}
    \min_{\mathbf{z}^{\rm{B}}} & \enskip g^{\text{obj}}(\mathbf{z}^{\rm{F}},\mathbf{z}^{\rm{B}}, \mathbf{y}) \label{eq:generic_realtime_obj} \\
    \text{s.t.} & \enskip g^{\text{gen}}(\mathbf{z}^{\rm{F}},\mathbf{z}^{\rm{B}}, \mathbf{y}) \leq 0 \label{eq:generic_realtime_gen} \\
    \phantom{s.t.} & \enskip g^{\text{net}}(\mathbf{z}^{\rm{F}},\mathbf{z}^{\rm{B}}, \mathbf{y}) \leq 0, \label{eq:generic_realtime_net}
\end{align}\label{eq:generic_realtime}
\end{subequations}
\noindent where $\mathbf{z}^{\rm{B}}$ symbolizes the adjustment variables in relation to the forward decisions $\mathbf{z}^{\rm{F}}$---here treated as a known and fixed input vector coming from~\eqref{eq:generic_forward}--- and $\mathbf{y}$ represents the eventually realized value of the uncertain parameters. The objective function~\eqref{eq:generic_realtime_obj} includes both the forward dispatch cost (which, at this point, is known and constant and, therefore, could be removed from the minimization) and the imbalance cost caused by unexpected net-load variations that occur very close to the real-time operation of the system. Similarly to problem~\eqref{eq:generic_forward}, the constraints~\eqref{eq:generic_realtime_gen} and \eqref{eq:generic_realtime_net} enforce the technical limitations of the generating units and the network, respectively. \textcolor{black}{As customary, the objective function and constraints of models \eqref{eq:generic_forward} and \eqref{eq:generic_realtime} are well known and properly defined.}

\cite{morales2014electricity} show that the value $\hat{\mathbf{y}}$ of the net demand that is used in the forward scheduling problem~\eqref{eq:generic_forward} may have a major impact on the subsequent balancing costs \eqref{eq:generic_realtime_obj}. Current practice, however, is content with a simple and direct solution, which is to take $\hat{\mathbf{y}}$ as an estimate of the expectation of the net demand $\mathbf{y}$ conditional on all the information at the forecaster's disposal. This information is generally known as \emph{context}. Yet, this expectation is oblivious to the minimization of the balancing costs that drives the real-time balancing problem~\eqref{eq:generic_realtime} and therefore, may turn out to be highly suboptimal. {\color{black} In other words, the conditional expectation overlooks the typical asymmetries affecting the minimization of the balancing costs, such as the fact that the distribution of the renewable generation is often skewed and that the costs of supplying upward and downward balancing energy are usually different.}


Instead of employing the conditional expectation of the uncertain parameters $\mathbf{y}$ as the estimate $\hat{\mathbf{y}}$ used in the forward scheduling problem~\eqref{eq:generic_forward}, we propose a regression procedure that provides an \emph{alternative value} for $\hat{\mathbf{y}}$ that explicitly accounts for the potential impact of the uncertain parameters on the subsequent real-time balancing problem~\eqref{eq:generic_realtime}. This alternative value is what is called a \emph{prescription} and a procedure to obtain it runs as follows. Suppose we have a sample of $N$ data points expressed in the form $\{(\mathbf{x}_i, \mathbf{y}_i)\}_{i \in \mathcal{N}} := \{(\mathbf{x}_1, \mathbf{y}_1), \ldots, (\mathbf{x}_i, \mathbf{y}_i), \ldots, (\mathbf{x}_N, \mathbf{y}_N)\}$, where $\mathbf{x} \in \mathbb{R}^p$ is a vector of features or covariates making up the \emph{context} and $\mathbf{y} \in \mathbb{R}^m$ is a vector of uncertain parameters. Our objective is to utilize said sample to infer a functional relation $\hat{\mathbf{y}} = h(\mathbf{x})$ with $h: \mathbb{R}^{p} \rightarrow \mathbb{R}^{m}$ such that, given the context $\mathbf{x}$, the provided prescription $\hat{\mathbf{y}}$ is trained to deliver the minimum total system costs in expectation when inserted into the forward scheduling problem~\eqref{eq:generic_forward}.

Now suppose that the prescriptive function $h$ is parameterized on a coefficient vector $\boldsymbol{q}$ of appropriate dimension, we propose to train $h_{\boldsymbol{q}}(\cdot)$, i.e., to estimate $\boldsymbol{q}$, by way of the following bilevel optimization problem:
\begin{subequations}
\begin{align}
    \min_{\boldsymbol{q}, \mathbf{z}_i^{\rm{F}}, \mathbf{z}_i^{\rm{B}}} & \enskip \frac{1}{|\mathcal{N}|}\sum_{i \in \mathcal{N}}  g^{\text{obj}} (\mathbf{z}_i^{\rm{F}}, \mathbf{z}_i^{\rm{B}}, \mathbf{y}_i) \label{eq:generic_bilevel_ul_obj} \\
    \text{s.t.} & \enskip g^{\text{gen}} (\mathbf{z}_i^{\rm{F}}, \mathbf{z}_i^{\rm{B}}, \mathbf{y}_i) \le 0, \enskip \forall i \in \mathcal{N} \label{eq:generic_bilevel_ul_gen} \\
    \phantom{s.t.} & \enskip g^{\text{net}} (\mathbf{z}_i^{\rm{F}}, \mathbf{z}_i^{\rm{B}}, \mathbf{y}_i) \le 0, \enskip \forall i \in \mathcal{N} \label{eq:generic_bilevel_ul_net} \\
    \phantom{s.t.} & \enskip \mathbf{z}_i^{\rm{F}} \in \{\argmin_{\mathbf{z}} \enskip f^{\text{obj}}(\mathbf{z}, h_{\boldsymbol{q}}(\mathbf{x}_i)) \label{eq:generic_bilevel_ll_obj} \\
    \phantom{s.t.} & \hspace{22mm} \text{s.t.}  \enskip f^{\text{gen}}(\mathbf{z}, h_{\boldsymbol{q}}(\mathbf{x}_i)) \le 0 \label{eq:generic_bilevel_ll_gen} \\
    \phantom{s.t.} & \hspace{22mm} \phantom{s.t.} \enskip  f^{\text{net}}(\mathbf{z}, h_{\boldsymbol{q}}(\mathbf{x}_i)) \le 0\}, \forall i \in \mathcal{N}. \label{eq:generic_bilevel_ll_net}
\end{align} \label{eq:generic_binivel}
\end{subequations}
Essentially, the bilevel program~\eqref{eq:generic_binivel} seeks the  $\boldsymbol{q}$-parameterized prescriptive function $h_{\boldsymbol{q}}(\cdot)$ that minimizes the empirical expectation of the total dispatch and re-dispatch costs over the data set $\{(\mathbf{x}_i, \mathbf{y}_i)\}_{i \in \mathcal{N}}$. For this purpose, it replicates, per data point $\{(\mathbf{x}_i, \mathbf{y}_i)\}$, the sequential two-stage scheduling process consisting of the forward power dispatch \eqref{eq:generic_bilevel_ll_obj}--\eqref{eq:generic_bilevel_ll_net} and the subsequent balancing re-dispatch~\eqref{eq:generic_bilevel_ul_obj}--\eqref{eq:generic_bilevel_ul_net}. {\color{black} From a statistical point of view, problem~\eqref{eq:generic_binivel} takes the form of an \emph{empirical risk minimization problem} and as such, is amenable to \emph{regularization} and \emph{robustification} \citep{shafieezadeh2019regularization}, admitting strategies for \emph{feature selection} too.}

Depending on the nature of variables $\mathbf{z}^{\rm{F}}, \mathbf{z}^{\rm{B}}$ and functions $f^{\text{obj}}$, $f^{\text{gen}}$, $f^{\text{net}}$, $g^{\text{obj}}$, $g^{\text{gen}}$, $g^{\text{net}}$, and $h_{\boldsymbol{q}}$, the difficulty of solving the bilevel optimization problem~\eqref{eq:generic_binivel} can vary significantly. For instance, if the forward scheduling decisions $\mathbf{z}^{\rm{F}}$ include binary variables representing the on/off status of  thermal generating units, then computing the global optimal solution of the bilevel problem is tremendously challenging as discussed in~\cite{fanghanel2009bilevel}. Likewise, even in the case that all scheduling decisions $\mathbf{z}^{\rm{F}}$ are continuous, the lower-level problem~\eqref{eq:generic_bilevel_ll_obj}--\eqref{eq:generic_bilevel_ll_net} can be non-convex if constraints $g^{\text{net}}$ account for the nonlinear AC power flow equations. In that case too, computing the global optimal solution of \eqref{eq:generic_binivel} is also a very challenging task. If all variables are continuous and all functions affine, then the bilevel optimization problem~\eqref{eq:generic_binivel} can be solved by replacing the lower-level problem with its equivalent KKT optimality conditions \citep{munoz2021bilevel}. Said strategy, however, requires additional binary variables and large enough constants that  increase the numerical instability and the computational burden of the resulting single-level optimization problem. \textcolor{black}{Besides, as discussed in \cite{pineda2019solving} and \cite{kleinert2020there}, existing methods to validate these upper bounds may fail and lead to suboptimal solutions.}

Given that the development of methods to solve a generic bilevel program like \eqref{eq:generic_binivel} is outside the scope of this paper, in the next section we  make use of some simplifying assumptions that allow us to efficiently solve the bilevel problem~\eqref{eq:generic_binivel} to global optimality using commercially available software for mixed-integer programming. {\color{black} In practice, we may not need to solve the estimation problem~\eqref{eq:generic_binivel} to global optimality. In effect, finding a good feasible solution to this problem may be enough to reap the bulk of the benefits of our approach. This enables the use of heuristics to solve problem~\eqref{eq:generic_binivel}.}

We also remark that some system operators are indeed aware of the importance of scheduling generation taking into account the subsequent real-time operation of the power system. For example, the California Independent System Operator (CAISO) \emph{manually modifies the load forecast}, anticipating the ramping capacity needs based on the operator's experience \citep{yurdakul2021forecasting}. This is explicitly explained in~\cite{CAISO} as follows: ``Operators in the ISO and energy imbalance market can manually modify load forecasts used in the market through load adjustments (...) to increase ramping capacity within the ISO by (...) committing additional units.'' Our approach is not only consistent with this \emph{modus operandi}, but also supports it with a scientifically grounded methodology to modify the load forecast.

Once the bilevel problem~\eqref{eq:generic_binivel} is solved and the prescriptive function $h_{\boldsymbol{q}}$ is obtained, the forward dispatch for an unseen value of the context $\mathbf{x}$ can be easily determined by solving the following deterministic optimization problem:
\begin{subequations}
\begin{align}
    \min_{\mathbf{z}^{\rm{F}}} & \enskip f^{\text{obj}}(\mathbf{z}^{\rm{F}},h_{\boldsymbol{q}}(\mathbf{x})) \label{eq:unseen_forward_obj} \\
    \text{s.t.} & \enskip f^{\text{gen}}(\mathbf{z}^{\rm{F}},h_{\boldsymbol{q}}(\mathbf{x})) \leq 0 \label{eq:unseen_forward_gen}\\
    \phantom{s.t.} & \enskip f^{\text{net}}(\mathbf{z}^{\rm{F}},h_{\boldsymbol{q}}(\mathbf{x})) \leq 0. \label{eq:unseen_forward_net}
\end{align}\label{eq:unseen_forward}
\end{subequations}
It is important to realize that the strategy we propose is completely in contrast with those other approaches that make use of stochastic programming to \emph{simultaneously} decide the forward scheduling decisions and the real-time adjustments, see, e.g., \cite{AGHAEI20091675, pritchard2010single, zavala2017stochastic}. Indeed, our approach builds upon the fact that, in practice, the forward dispatch and the subsequent balancing re-dispatch are \emph{sequentially} optimized. The next section elaborates on the simplifying hypotheses we make to transform the bilevel program  \eqref{eq:generic_binivel} into a single-level mixed-integer optimization problem.

\section{Contextual merit-order dispatch}\label{sec:CD}
For the remainder of this article, we assume the following.

\begin{itemize}
    \item Forward dispatch decisions and real-time balancing actions are all modeled by continuous variables, that is, on/off binary variables are not considered. 
    \item The inter-temporal constraints of power production portfolios, such as ramping limits and minimum-up and -down times, are not explicitly accounted for in the generation scheduling problem.
    \item Network constraints are only taken care of in the real-time stage using a pipeline representation of the transmission network. That is, constraints~\eqref{eq:generic_bilevel_ll_net} are removed from the bilevel optimization problem~\eqref{eq:generic_binivel}.
    \item All power loads are assumed to be non-dispatchable and therefore, the objective function of the two-stage scheduling problem boils down to the minimization of the generating cost of dispatchable units. 
\end{itemize}

With these assumptions in place, the forward scheduling problem~\eqref{eq:generic_forward} is simplified to a merit-order dispatch whereby generators are scheduled based on their marginal production costs, as explained later. Furthermore, {\color{black} with the pipeline representation of the transmission network,} our setup becomes more aligned with European practices, in keeping with the case study we present in Section~\ref{sec:CaseStudy}. Likewise, we consider non-dispatchable loads only just to simplify the subsequent exposition of our approach. Finally, by way of the rest of the assumptions, the bilevel model~\eqref{eq:generic_binivel} becomes simpler and computationally manageable.  

Nevertheless, any of these assumptions could {\color{black} also} be dropped at the expense of increasing the complexity of this model as discussed at the end of Section \ref{sec:PD}. {\color{black} For instance, if on/off commitment decisions of the generating units are considered, solving the bilevel model~\eqref{eq:generic_binivel} requires specific solution methodologies depending on the application. If all decisions are continuous but additional intertemporal or network constraints are included in the day-ahead stage, then the single-level reformulation of model~\eqref{eq:generic_binivel} requires dual variables of the lower-level constraints and large enough constants, which increases its computational burden and endangers optimality guarantees. Finally, the use of a linearized (e.g., DC) power flow model in the real-time operational stage~\eqref{eq:generic_realtime} of the generation scheduling problem would not have any qualitative impact on the computational complexity of the estimation problem~\eqref{eq:generic_binivel}.}

All in all, the particular forward generation scheduling problem we consider in this paper is given as follows:
\begin{subequations}
\begin{align}
    \min_{p_g, g \in G} & \enskip \sum_{g \in G} C_g p_{g} \\
    \text{s.t.} & \enskip \sum_{g \in G}  p_{g} = \widehat{L} \label{eq:merit_order_eq}\\
    \phantom{s.t.} & \enskip 0 \le p_{g} \le \overline{P}_{g}, \enskip \forall g \in G, \label{eq:merit_order_block_limit}
\end{align}\label{eq:merit_order}
\end{subequations}
%
where $p_g, \overline{P}_{g}, C_g  \in \mathbb{R}^{+}$  and $G \subseteq \mathbb{N}$ is the set of generation units (it can also represent generation blocks with different cost). Each unit $g$ has associated a production level $p_g$ and a marginal cost $C_g$. Equation~\eqref{eq:merit_order_eq} enforces the aggregate power balance, with parameter $\widehat{L} \in \mathbb{R^+}$ representing a point or single-value estimate of the \emph{total net demand} $L\in \mathbb{R^+}$ in the system, which is unknown at the moment the forward scheduling is determined and thus, is to be treated as a random variable. 
Lastly, equation~\eqref{eq:merit_order_block_limit} sets the capacity of each generation unit.

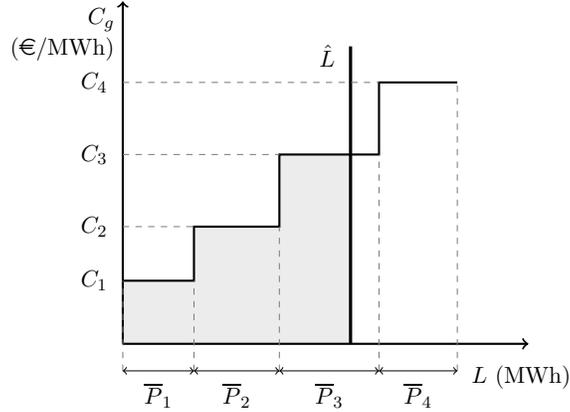
\begin{figure} 
\def\pga{1.0} 
\def\pgb{2.2}
\def\pgc{3.6}
\def\pgd{4.7}
\def\ca{0.7} 
\def\cb{1.3}
\def\cc{2.1}
\def\cd{2.9}
\def\psum{3.2}
\def\arrowy{-0.3}
\centering
    \begin{tikzpicture}[scale=0.8]
    	\begin{axis}[
    	    scale=1,
        	hide axis,
        	width=14cm,
            xmin = -2.5, xmax = 8, ymin=-1, ymax=6,
        	clip marker paths=true]

            \addplot+[mark=none, draw=none, fill=lightgray!30] coordinates {(0, \ca) (\pga, \ca)} \closedcycle;
            \addplot+[mark=none, draw=none, fill=lightgray!30] coordinates {(\pga, \cb) (\pgb, \cb)} \closedcycle;
            \addplot+[mark=none, draw=none, fill=lightgray!30] coordinates {(\pgb, \cc) (\psum, \cc)} \closedcycle;
            
        	\addplot[mark=none, solid, line width=1pt, ->] coordinates {(0, 0) (5.7,0)};
        	\addplot[mark=none, solid, line width=1pt, ->] coordinates {(0, 0) (0, 3.8)};
        	\node[anchor=east] (pp) at (axis cs:0, 3.6){$C_g$};
        	\node[anchor=east] (ppb) at (axis cs:0, 3.27){(\euro/MWh)};
            \node[anchor=north] (qq) at (axis cs:5.6, -0.15){$L$ (MWh)};
            
        	\addplot[mark=none, solid, domain=0:5, line width=1pt] coordinates { (0, \ca) (\pga, \ca) (\pga, \cb) (\pgb, \cb) (\pgb, \cc) (\pgc, \cc) (\pgc, \cd) (\pgd, \cd) };
        	
        	\addplot[mark=none, solid, <->] coordinates {(0,\arrowy) (\pga,\arrowy)};
        	\node[anchor=north] (db) at (axis cs: 0.5 , \arrowy - 0.05){$\overline{P}_1$};
        	\addplot[mark=none, solid, <->] coordinates {(\pga,\arrowy) (\pgb,\arrowy)};
        	\node[anchor=north] (db) at (axis cs: 1.6 , \arrowy - 0.05){$\overline{P}_2$};
        	\addplot[mark=none, solid, <->] coordinates {(\pgb,\arrowy) (\pgc,\arrowy)};
        	\node[anchor=north] (db) at (axis cs: 2.9 , \arrowy - 0.05){$\overline{P}_3$};
        	\addplot[mark=none, solid, <->] coordinates {(\pgc,\arrowy) (\pgd,\arrowy)};
        	\node[anchor=north] (db) at (axis cs: 4.15 , \arrowy - 0.05){$\overline{P}_4$};
            
            
            \addplot[mark=none, dashed, color=gray, line width=0.5pt] coordinates {(0, \arrowy) (0,\ca)};
            \addplot[mark=none, dashed, color=gray, line width=0.5pt] coordinates {(\pga, \arrowy) (\pga,\ca)};
            \addplot[mark=none, dashed, color=gray, line width=0.5pt] coordinates {(\pgb, \arrowy) (\pgb,\cb)};
            \addplot[mark=none, dashed, color=gray, line width=0.5pt] coordinates {(\pgc, \arrowy) (\pgc,\cc)};
            \addplot[mark=none, dashed, color=gray, line width=0.5pt] coordinates {(\pgd, \arrowy) (\pgd,\cd)};
            \addplot[mark=none, solid, line width=1.5pt] coordinates {(\psum, 0) (\psum,3.3)};
            \node[anchor=east] (xx) at (axis cs: \psum - 0.1 , \cd + 0.3){$\hat{L}$};

            \addplot[mark=none, dashed, color=gray, line width=0.5pt] coordinates {(0, \cb) (\pga,\cb)};
			\addplot[mark=none, dashed, color=gray, line width=0.5pt] coordinates {(0, \cc) (\pgb,\cc)};
			\addplot[mark=none, dashed, color=gray, line width=0.5pt] coordinates {(0, \cd) (\pgc,\cd)};

            \node[anchor=east] (ccc) at (axis cs: -0.1 , \ca){$C_1$};
            \node[anchor=east] (ccc) at (axis cs: -0.1 , \cb){$C_2$};
            \node[anchor=east] (ccc) at (axis cs: -0.1 , \cc){$C_3$};
            \node[anchor=east] (ccc) at (axis cs: -0.1 , \cd){$C_4$};
            
    	\end{axis}	
    \end{tikzpicture} \\
\caption{{\color{black} Illustration of a \emph{cost-merit order} dispatch. The vertical solid line represents the total net demand $\widehat{L}$ and the shadow area represents the optimal dispatch computed by model~\eqref{eq:merit_order}}.}
\label{fig:merit_order}
\end{figure}

The linear program~\eqref{eq:merit_order} stands for an economic dispatch problem whereby generation units are dispatched following a cost-merit order, meaning that units $g$ with a lower cost $C_g$ are dispatched first. \textcolor{black}{This is illustrated in Figure~\ref{fig:merit_order} for four units with capacities $\overline{P}_1, \overline{P}_2, \overline{P}_3, \overline{P}_4$ and marginal costs $C_1, C_2, C_3, C_4$, respectively. The vertical solid line represents the total net demand $\widehat{L}$ and the shadow area indicates the optimal dispatch computed by model~\eqref{eq:merit_order}. Since units 1 and 2 are the cheapest ones, their dispatch is set at their maximum capacity. Unit 3 is partially dispatched until the demand is fulfilled. Unit 4 is the most expensive one and therefore, is left out of the demand supply.} To ease the discussion that follows, hereinafter we consider that the units in the set $G$ are ordered such that $g < g'$ if and only if $C_g < C_{g'}$. Hence, if we denote the optimal solution to~\eqref{eq:merit_order} as $\{p^*_g\}_{g \in G}$, it holds that $p^*_{g'} >0$ implies that $p^*_{g} = \overline{P}_{g}$, whenever $g < g'$.

Since the system \emph{net} demand is uncertain, the power dispatch that results from the forward problem~\eqref{eq:merit_order} is to be adjusted during the real-time operation of the power system to satisfy the \emph{actual} net demand. The aim of the real-time balancing problem is to correct the imbalance of the system in a cost-efficient manner. To this end, we consider a pipeline model where $G(b)$ and $D(b)$ represent the set of generation units and loads that are connected to node $b$, in that order. With some abuse of notation, let $(L_{di}\in \mathbb{R})_{d \in D(b)}$ be a certain realization $i \in \mathcal{N}$ of the net load $d \in D(b)$ connected to node $b$ of the system. We also define $o(l)$ and $e(l)$ as the origin and ending nodes of line $l$, respectively. Thus, $\{l: o(l) = b\}$ and $\{l: e(l) = b\}$ represent the subset of lines that start or end at node $b$, in that order. Once introduced this notation, the real-time balancing problem under consideration renders
\begin{subequations}
\begin{align}
    \min_{\Xi} & \enskip \sum_{g=1}^{G}(C_g^{\text{u}}r_{g}^{\text{u}} - C_g^{\text{d}} r_{g}^{\text{d}})\label{eq:real_time_OF} \\
    \text{s.t.} & \enskip 0\le p^{*}_{g} + r_{g}^{\text{u}} - r_{g}^{\text{d}} \le \overline{P}_{g}, \enskip \forall g \in G \label{eq:real_time_BlockSize} \\
    \phantom{s.t.} & \enskip 0 \le r_{g}^{\text{u}} \le R_g^{\text{u}}, \enskip \forall g \in G\label{eq:real_time_URsize} \\
    \phantom{s.t.} & \enskip 0 \le r_{g}^{\text{d}} \le R_g^{\text{d}}, \enskip \forall g \in G\label{eq:real_time_DRsize} \\
    \phantom{s.t.} & \enskip \sum_{g \in G(b)}(p^{*}_{g} + r_{g}^{\text{u}} - r_{g}^{\text{d}}) = \sum_{d \in D(b)} L_{di} + \sum_{l:o(l)=b} f_{l} - \sum_{l:e(l)=b} f_{l}, \enskip \forall b \in B\label{eq:real_time_nodalPB} \\
    \phantom{s.t.} & \enskip |f_{l}| \le \overline{F}_{l}\label{eq:real_time_LineCapacity}, \enskip \forall l \in \Lambda,
\end{align} \label{eq:real_time}
\end{subequations}
where $\Xi:=\{r_{g}^{\text{u}}, r_{g}^{\text{d}}  \in \mathbb{R}^{+}, g \in G, f_l  \in \mathbb{R}, l \in \Lambda\}$ is the set of decision variables and $p^{*}_{g}, R_g^{\text{u}}, R_g^{\text{d}}, C_g^{\text{u}}, C_g^{\text{d}}, \overline{P}_{g}, \overline{F}_{l} \in  \mathbb{R}^{+}$ and $L_{di} \in \mathbb{R}$ are known parameters.

The power output of each flexible unit $g$ may be increased by an amount $r_{g}^{\rm u}$, based on the marginal cost for upward balancing energy $C_g^{\text{u}}$, or decreased by an amount $r_{g}^{\rm d}$ of downward balancing energy, which entails a marginal benefit (linked to fuel-cost savings) of $C_g^{\text{d}}$. These actions are driven by the nodal power balance equation~\eqref{eq:real_time_nodalPB} and the minimization of the total balancing costs~\eqref{eq:real_time_OF}. Naturally, the amount of balancing energy provided from each generation unit $g$, either upward or downward, must be such that the eventual power output from that unit (taking into account the forward optimal schedule $p_{g}^{*}$) is positive and lower than its capacity $\overline{P}_{g}$, as stated in equation~\eqref{eq:real_time_BlockSize}.  Moreover, constraints~\eqref{eq:real_time_URsize} and~\eqref{eq:real_time_DRsize} limit the amount of up- and down-balancing power that can be deployed from each generation unit to $R_g^{\text{u}}$ and $R_g^{\text{d}}$, which are indicative of how flexible the underlying asset actually is. Finally, line capacity limits are imposed by~\eqref{eq:real_time_LineCapacity}, with $f_l$ being the power flow through line $l$ and $\overline{F}_l$ the capacity of the line.


Now suppose we have a sample of $N$ data points expressed in the form $\{(\mathbf{x}_i, L_i)\}_{i \in \mathcal{N}} := \{(\mathbf{x}_1, L_1), \ldots, (\mathbf{x}_i, L_i), \ldots, (\mathbf{x}_N, L_N)\}$, where $\mathbf{x} \in \mathbb{R}^p$ is a vector of features or covariates making up the \emph{context} and $L \in \mathbb{R}^+$ is the random net system demand. As stated in the previous section, our objective is to utilize said sample to infer a functional relation $\widehat{L} = h(\mathbf{x})$, $h: \mathbb{R}^{p} \rightarrow \mathbb{R}^+$, such that, given the context $\mathbf{x}$, the provided prescription $\widehat{L}$ is trained to deliver the minimum total system costs in expectation when inserted into the power balance equation~\eqref{eq:merit_order_eq}. For simplicity, and because it proves to perform very satisfactorily in the numerical experiments of Section~\ref{sec:CaseStudy}, we restrict $h$ to the family of affine linear functions, i.e., $h(\mathbf{x}) = \mathbf{q}^{\top}\, \mathbf{x}$, with $\mathbf{q} \in \mathbb{R}^{p}$ and with one of the features, say $x_1$, fixed to one. This selection, together with the particular structure of problems~\eqref{eq:merit_order} and \eqref{eq:real_time} circumvent the need for applying the conventional procedures to solve the bilevel problem~\eqref{eq:generic_binivel}. Instead, to determine $\mathbf{q}$, we solve the following empirical risk minimization problem:
%
\begin{subequations}
\begin{align}
    \min_{\mathbf{q},\, \Upsilon} & \enskip \frac{1}{N}\sum_{i \in \mathcal{N}} \sum_{g \in G} (C_g p_{gi} + C_g^{\text{u}} r_{gi}^{\text{u}} - C_g^{\text{d}} r_{gi}^{\text{d}}) \label{eq:estimation_problem_of}\\
    \text{s.t.} & \enskip 0 \le p_{gi} + r_{gi}^{\text{u}} - r_{gi}^{\text{d}} \le \overline{P}_{gi}, \enskip \forall i \in \mathcal{N}, \enskip \forall g \in G \label{eq:estimation_problem_two_stage_balance} \\
    \phantom{s.t.} & \enskip 0 \le r_{gi}^{\text{u}} \le R_g^{\text{u}}, \enskip \forall i \in \mathcal{N}, \enskip \forall g \in G \\
    \phantom{s.t.} & \enskip 0 \le r_{gi}^{\text{d}} \le R_g^{\text{d}}, \enskip \forall i \in \mathcal{N}, \enskip \forall g \in G\\
    \phantom{s.t.} & \enskip \sum_{g \in G(b)}(p_{gi} + r_{gi}^{\text{u}} - r_{gi}^{\text{d}}) = \hspace{-2mm} \sum_{d \in D(b)} \hspace{-2mm} L_{di} + \hspace{-2mm}\sum_{l:o(l)=b} \hspace{-2mm} f_{li} - \hspace{-2mm} \sum_{l:e(l)=b} \hspace{-2mm} f_{li}, \forall i \in \mathcal{N}, \forall b \in B\\
    \phantom{s.t.} & \enskip |f_{li}| \le \overline{F}_{l}, \enskip \forall i \in \mathcal{N}, \enskip \forall l \in \Lambda \label{eq:estimation_problem_linecapacity}\\
    \phantom{s.t.} & \enskip \sum_{g \in G}  p_{gi} =  \widehat{L}_{i}, \enskip \forall i \in \mathcal{N} \label{eq:estimation_problem_PB1}\\
    \phantom{s.t.} & \enskip \widehat{L}_{i} = \sum_{j=1}^{p} q_{j} x_{ji}, \enskip \forall i \in \mathcal{N} \label{eq:estimation_problem_AF}\\
    \phantom{s.t.} & \enskip u_{gi} \overline{P}_g \le p_{gi} \le u_{(g-1)i}\, \overline{P}_g, \enskip \forall i \in \mathcal{N}, \enskip \forall g \in G: g > 1 \label{eq:estimation_problem_MO1}\\
    \phantom{s.t.} & \enskip u_{gi} \overline{P}_g \le p_{gi} \le \overline{P}_g, \enskip \forall i \in \mathcal{N}, \enskip g = 1\\
    \phantom{s.t.} & \enskip u_{gi} \le u_{(g - 1)i}, \enskip \forall i \in \mathcal{N}, \enskip \forall g \in G: g > 1 \\
    \phantom{s.t.} & \enskip u_{gi} \in \{0, 1\}, \enskip \forall i \in \mathcal{N}, \enskip \forall g \in G, \label{eq:estimation_problem_MOfinal}
\end{align} \label{eq:estimation_problem}
\end{subequations}
where $\Upsilon:=\{p_{gi}, r_{gi}^{\text{u}}, r_{gi}^{\text{d}}, u_{gi}, f_{li}\}_{\{i, g,l\}}$. Intuitively, the estimation problem~\eqref{eq:estimation_problem} computes the coefficient vector $\mathbf{q} = (q_{1}, \ldots, q_{p})$ such that the total system cost averaged over the sample is minimized. This is the reason why all the decision variables related to the power dispatch and the provision of regulating power, i.e., $p_{gi}, r_{gi}^{\text{u}}, r_{gi}^{\text{d}}$, appear augmented with the sample index $i$ in~\eqref{eq:estimation_problem}. Constraints~\eqref{eq:estimation_problem_two_stage_balance}--\eqref{eq:estimation_problem_PB1} serve exactly the same purpose as their analogs in \eqref{eq:merit_order} and \eqref{eq:real_time}. Equation~\eqref{eq:estimation_problem_AF} expresses the prescription $\widehat{L}_{i}$ of the net system demand $L$ under context $\mathbf{x}_i$ as an affine function of the features, whose coefficients are to be computed by solving~\eqref{eq:estimation_problem}. Finally, the set of constraints~\eqref{eq:estimation_problem_MO1}--\eqref{eq:estimation_problem_MOfinal} guarantee that the power dispatch $\{p_{gi}\}_{g\in G}$ coming from~\eqref{eq:estimation_problem} for each sample $i$ is optimal in the forward problem~\eqref{eq:merit_order}, that is, these constraints constitute the \emph{optimality conditions} of the forward dispatch problem~\eqref{eq:merit_order}. Accordingly, these constraints enforce, for each sample $i$, the merit-order dispatch of the generation units $\{p_{gi}\}_{g\in G}$, forcing that $p_{g'i} >0 \implies p_{gi} = \overline{P}_{g}$, for all $g \in G: g < g'$. Importantly, this formulation of the optimality conditions of problem~\eqref{eq:merit_order} neither necessitates dual variables, nor large enough constants, both of which are frequently used to solve lineal bilevel programs \citep{pineda2019solving}. 

Problem~\eqref{eq:estimation_problem} is, therefore, a mixed-integer linear program due to the binary character of variables $u_{gi}$, which are used to impose the cost-merit order. As such, this problem can be solved using commercially available solvers such as CPLEX \citep{CPLEX}. Once we obtain the optimal coefficient vector $\mathbf{q}^*$, we can produce the net-demand prescription $\widehat{L} = (\mathbf{q}^*)^{\top}\mathbf{x}$, which is to be fed into~\eqref{eq:merit_order_eq} under the context $\mathbf{x}$ to readily obtain the forward dispatch decisions.

Now, let $L^{\rm F}$ denote the expected net-demand value (which is typically known as a \emph{point prediction} of the net demand). Even though this expected value has been and is normally used as $\widehat{L}$ in the forward problem~\eqref{eq:merit_order}, it is \emph{not} consistent with the plausible asymmetry in the cost of dealing with the subsequent prediction errors through the real-time problem~\eqref{eq:real_time}. Indeed, it is most often the case that the cost of increasing the electricity production in real time is different from that of diminishing it. In this line, problem~\eqref{eq:estimation_problem} offers a handy way to construct a new estimate $\widehat{L}$ to be used in \eqref{eq:merit_order} that takes into account the referred cost asymmetry. Indeed, the training problem~\eqref{eq:estimation_problem} can be used in practice to upgrade the point prediction $L^{\rm F}$ to a prescription for $\widehat{L}$, which does account for the asymmetry of the power system's balancing costs. For this, it suffices to include $L^{\rm F}$ as one of the features that are part of the context $\mathbf{x}$. This is what we do in the following example and in the case study of Section~\ref{sec:CaseStudy}. {\color{black} That said, the context $\boldsymbol{x}$ could include any other information that could be deemed useful to improve the prescription $\widehat{L}$, for instance, the estimated standard deviation or quantiles of the conditional net demand, if these were available through a probabilistic forecast.}

We finish this discussion with an important remark. Despite the fact that constraints~\eqref{eq:estimation_problem_MO1}--\eqref{eq:estimation_problem_MOfinal} turn our training problem~\eqref{eq:estimation_problem} into a mixed-integer program,  enforcing the cost-merit order through these constraints is critical to train an affine model $h(\mathbf{x}) = \mathbf{q}^{\top}\mathbf{x}$ that renders economic benefits within the two-stage scheduling problem described in Section~\ref{sec:PD}. Furthermore, even though the training problem~\eqref{eq:estimation_problem} may require some computational effort to be solved, the affine model it delivers is intended to remain effective for a period of time (e.g., weeks or months), and hence, the task of solving the mixed-integer program~\eqref{eq:estimation_problem} only has to be undertaken once in a while. {\color{black} That said, the parameter vector $\mathbf{q}$ is to be re-estimated \emph{offline} on a regular basis to capture changes in the dynamics of the electricity market and the power system (for example, to account for seasonal variations).}

\section{Example}\label{sec:Ex}
\begin{figure}[htp]
    \centering
    \begin{tikzpicture}
    \draw
    (0,0) node [sin v source] (v1) {} 
    (v1.north)--++(0,0.5) coordinate(v1-gen) 
    ($(v1-gen) + (0, -1.25) $) node[below]{$G_1$}
    (v1-gen)--++(-0.5,0) 
    (v1-gen)--++(1.5,0) node[right]{1} 
    coordinate[pos=0.5](v1-r) 
    (v1-r)--++(0,-0.25)--++(2.625,-2)--++(0,-0.25) 
    coordinate(v3-l) 
    ($(v1-r) + (2, -1) $) node{Line 1}
    (v3-l)--++(-0.25,0) node[left]{3} 
    (v3-l)--++(1,0) 
    coordinate[pos=0.375](v3-dem) 
    ($(v3-dem) + (0, -0.6) $) node[below]{$L$}
    coordinate[pos=0.75](v3-r) 
    (v3-r) --++(0,0.25)--++(2.625,2)--++(0,0.25) 
    coordinate(v2-l)
    (v2-l)--++(-0.75,0) node[left]{2} 
    (v2-l)--++(1.25,0) 
    coordinate[pos=0.6](v2-gen) 
    (v2-gen) --++(0,-0.5) node[below,sin v source](v2){} 
    ($(v2-gen) + (0, -1.25) $) node[below]{$G_2$}
    ($(v2-l) + (-2, -1) $) node{Line 2}
    ;
    \draw[-stealth, thick](v3-dem)--++(0,-0.5cm); 
    \draw[ultra thick]  ($(v3-dem) + (-1.25 *.5, 0) $) -- ($(v3-dem) + (1.25 *.5, 0) $);
    \draw[ultra thick]  ($(v1-gen) + (-.5, 0) $) -- ($(v1-gen) + (1.5, 0) $);
    \draw[ultra thick]  ($(v2-gen) + (.5, 0) $) -- ($(v2-gen) + (-1.5, 0) $);
    \end{tikzpicture} 
    \caption{Three-bus power system with one random demand and two thermal generators. }
    \label{fig:toy_system}
\end{figure}
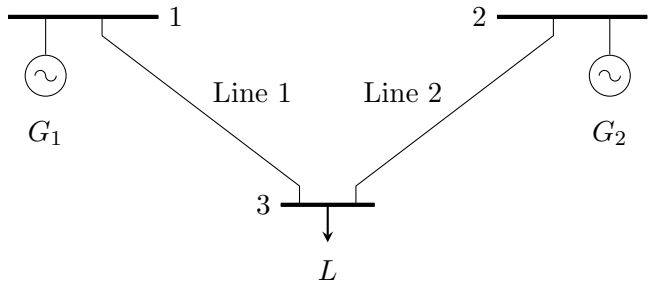

Consider the small three-bus system depicted in Figure~\ref{fig:toy_system}, which is composed of one random demand $L$ at bus 3, two thermal generators, $G_1$ and $G_2$, at buses 1 and 2, respectively; and two lines, Line 1 and Line 2, connecting nodes 1 and 3, and buses 2 and 3, in that order.

The technical and economic characteristics of generating units $G_1$ and $G_2$ are collated in Table~\ref{table:toy_gen_data}. Note that, in comparison, unit $G_1$ is smaller and cheaper than $G_2$. In contrast, the latter is significantly more flexible as it features re-dispatch costs, i.e., $C_g^{\text{u}}$ and $C_g^{\text{d}}$, that are much more competitive. We remark that $C_{1}^{\text{d}} =$ $-20$ \euro/MWh implies that this power unit must be paid 20 \euro \, for each MWh its production is decreased in the real-time stage. Unless stated otherwise, line capacities are assumed infinite.

The only demand in the system, namely, $L$, is random. Suppose we have a sample $\{(\mathbf{x}_i,L_{i})\}_{i=1}^{N}$, where the feature $\mathbf{x}_i = (1, L^{\rm F}_i)^{\top}$. Again, $L^{\rm F}_i$ represents a classical point prediction of the demand $L$ built out of whichever available information the forecaster had at her disposal to produce it by way of whatever machine learning or forecasting technique she could have developed to that end. We stress that this setup is very common in reality, where power system operators often count on specialized software to produce good point predictions $L^{\rm F}_i$. Our objective is to use our methodology and training model~\eqref{eq:estimation_problem} described in Section~\ref{sec:CD} to recycle this standard point prediction with the aim of fabricating a better value for $\widehat{L}$ in equation~\eqref{eq:merit_order_eq}.
\begin{table}
    \setlength{\tabcolsep}{12pt}
    \centering
    \begin{tabular}{cccccccc}
    \toprule
       & $C_g$           & $C_g^{\text{u}}$       & $C_g^{\text{d}}$       &  $\overline{P}_g$     & $R_g^{\text{u}}$  & $R_g^{\text{d}}$ \\
    \midrule
    $G_1$   & 5   & 30  & -20  &  60  &  60  &  60  \\
    $G_2$   & 15  & 20  &  10  & 150  & 150  & 150  \\
    \bottomrule
    \end{tabular}
    \captionsetup{justification=centering}
    \caption{Technical and economic specifications of power plants. Marginal costs are given in \euro/MWh and capacities in MW.}
    \label{table:toy_gen_data}
\end{table}

For this small example, we generate samples in the form $\{(\mathbf{x}_i, L_{i})\}_{i=1}^{N}$ as follows. We consider that the per-unit (p.u.) point forecast of the net demand $L$ follows a uniform distribution between $a$ and $b$. Therefore, $L^{\rm F} \sim \overline{L} \cdot U(a,b)$, where $\overline{L}$ is a factor representing the maximum power load at bus~3. We further assume that the per-unit net demand itself $L/\overline{L}$  follows a Beta distribution with mean equal to $L^{\rm F}/\overline{L}$ and standard deviation $\sigma$. Hence, $L \sim \overline{L}\cdot {\rm Beta}(\alpha, \beta)$, where the scale and shape parameters $\alpha$ and $\beta$ are related to the mean $L^{\rm F}/\overline{L}$ and the standard deviation $\sigma$ as follows:
\begin{subequations}
\begin{align}
    & \frac{L^{\rm F}}{\overline{L}} = \frac{\alpha}{\alpha + \beta}, \\
    & \sigma^2 = \frac{\alpha \cdot \beta}{(\alpha + \beta)^2 (\alpha + \beta + 1 )}.
\end{align} \label{eq:alpha_beta_1}
\end{subequations}
{\color{black} Equations~\eqref{eq:alpha_beta_1} guarantee that the point prediction $L^{\rm F}$ is an unbiased estimator of the conditional net demand $L$.}
In this illustrative example we fix $\sigma =0.075$ p.u. and generate 20 samples $\{(\mathbf{x}_i,L_{i})\}_{i=1}^{N}$ with $N =750$. Each $L^{\rm F}_i$ in $\mathbf{x}_i$ is randomly drawn from $\overline{L} \cdot U(a,b)$. Given $L_i^{\rm F}/\overline{L}$ and $\sigma$, and provided that the system of nonlinear equations~\eqref{eq:alpha_beta_1} has a solution (notice that $\alpha,\beta >0$), parameters $\alpha_i$ and $\beta_i$ can be computed as
\begin{subequations}
\begin{align}
    & \alpha_i = - \frac{1}{\sigma^2} \left(\left(\frac{L_i^{\rm F}}{\overline{L}}\right)^2 - \frac{L_i^{\rm F}}{\overline{L}} + \sigma^2\right)  \frac{L_i^{\rm F}}{\overline{L}}, \\
    & \beta_i = \frac{1}{\sigma^2} \left(\left(\frac{L_i^{\rm F}}{\overline{L}}\right)^2 - \frac{L_i^{\rm F}}{\overline{L}} + \sigma^2\right) \left(\frac{L_i^{\rm F}}{\overline{L}} - 1\right).
\end{align} \label{eq:alpha_beta_2}
\end{subequations}
Each $L_{i}$ is then randomly taken from $\overline{L}\cdot {\rm Beta}(\alpha_i,\beta_i)$. We take the first 500 data points of each sample as the training set and the last 250 as the test set.
%

We postulate the affine model $\widehat{L} = q_0 + q_1  L^{\rm F}$  and solve problem~\eqref{eq:estimation_problem} on the training set to determine coefficients $q_0$ and $q_1$. Finally, to evaluate the performance of the affine model, for each data point $(\mathbf{x}_i, L_{i})$ in the test set, we simulate the \emph{sequential} forward and real-time problems~\eqref{eq:merit_order} and~\eqref{eq:real_time}, with $\widehat{L} = q_0 + q_1 L^{\rm F}_i$ in~\eqref{eq:merit_order_eq}, and $L_{i}$ in~\eqref{eq:real_time_nodalPB}. We then compute the sum of the forward and real-time production costs averaged over the 250 data points in the test set. This mean sum is further averaged over the 20 samples we generate. Our approach, which uses a \textit{prescription} of the system net demand for scheduling generation, is referred to as P-SC (from prescriptive scheduling). We compare it with the customary practice of directly using the point \textit{forecast}  $ L^{\rm F}_i$ as $\widehat{L}$ in~\eqref{eq:merit_order_eq}, which is referred to as F-SC (from forecast scheduling). Notice that our approach boils down to the conventional one if $q_0=0$ and $q_1=1$. Finally, the relative cost difference between these approaches is denoted as $\Delta$cost. 

In the results we discuss next, we set a base case with $a=0.03$, $b=0.97$, $\overline{L} = 100$ MW, and the technical and economic parameters of the three-bus system described above\footnote{We take $a = 0.03$ and $b=0.97$ pu to ensure that \eqref{eq:alpha_beta_1} has a real solution.}. We then define variants of this case by changing one or some of those parameters.

\subsection{Impact of power regulation costs}\label{IE:reg_cost}
Table~\ref{table:toy_exp_1} provides the cost savings that our approach achieves with respect to the conventional one under different $G_2$'s power regulation costs. For completeness, this table also includes the average cost of these two approaches for the test set and the values of the intercept $q_0$ and the linear coefficient $q_1$ of the affine model for $\widehat{L}$ our approach utilizes. These values represent expectations over the test data points of the 20 samples generated as indicated above. The first row in the table corresponds to the base case.

Interestingly, our approach systematically corrects the point forecast of the net demand $L$ \emph{downwards}, with a linear coefficient $q_1$ which is, on average, lower than or equal to 1, and a negative intercept $q_0$ in expectation. This is so because it is economically advantageous for the system to cope with positive net demand errors (i.e., eventual demand increases) by deploying upward balancing energy from unit $G_2$. Indeed, the alternative would be to deal with negative demand errors by down-regulating (i.e., by providing downward balancing energy) with unit $G_1$, a recourse that is clearly much more expensive. 

To further elaborate on this phenomenon, the second row in Table~\ref{table:toy_exp_1} provides results for a variant of the base case in which $C_2^{u}$ has been decreased from 20 to 15 \euro/MWh. Now that up-regulating (i.e., providing upward balancing energy) through unit $G_2$ is even cheaper, the downward correction of our approach to the net demand point forecast is more pronounced and the associated cost savings due to said phenomenon become larger.  On the contrary, if it is the provision of downward balancing energy by $G_2$ what becomes 5 \euro/MWh cheaper and, hence, free (see third row of Table~\ref{table:toy_exp_1}), the net demand point forecast is barely corrected and the costs savings brought by our approach (with regard to F-SC) become smaller as a result. Note that correcting the point forecast \emph{upwards} in this case (in an attempt to profit from the free downward balancing energy provided by $G_2$) would be counterproductive in reality, as the system may risk having to resort to the high-cost downward balancing energy of unit $G_1$ in those likely scenarios in which the net demand ends up being lower than the capacity of this unit.
%

 At this point, it may be instructive to see what happens when we drop constraints \eqref{eq:estimation_problem_MO1}--\eqref{eq:estimation_problem_MOfinal} from the mixed-integer program through which we train the affine model $\widehat{L} = q_0 + q_1  L^{\rm F}$. This is indeed very tempting, because, if these constraints are removed, the training model~\eqref{eq:estimation_problem} becomes a very pleasant linear program, similar to the stochastic-programming-based formulation advocated, for instance, in~\cite{pritchard2010single}(with the data points in~\eqref{eq:estimation_problem} playing the role of the ``scenarios'' in~\cite{pritchard2010single}). However, these constraints ensure the optimality of the lower-level problem \eqref{eq:generic_bilevel_ll_obj}--\eqref{eq:generic_bilevel_ll_net} and guarantee that the above affine model is learned following a cost-merit-order principle.
 Therefore, if these constraints are dropped from \eqref{eq:estimation_problem}, the affine model $\widehat{L} = q_0 + q_1  L^{\rm F}$ is not trained for the target task. This is exactly what Table~\ref{table:toy_exp_1_no_merit} shows. This table is analogous to Table~\ref{table:toy_exp_1}, but for a linear training model made up of constraints \eqref{eq:estimation_problem_of}--\eqref{eq:estimation_problem_AF} only. We denote this approach as L-SC from ``Linear''). The training model L-SC ignores the merit order and hence,  takes for granted that the system can benefit from the cheap downward balancing energy of unit $G_2$ by allocating a non-zero production to this unit in the forward problem regardless of whether unit $G_1$ has been fully dispatched or not. This is, however, an strategy forbidden by the problem, which explains the poor actual performance of L-SC. This phenomenon is especially notorious  for the case $C_2^{d} = 15$ \euro/MWh, in which the demand is heavily overestimated (the mean of the estimated demand is increased by 45\%) and only the free downward regulation from unit $G_2$ is used.

\begin{table}
    \setlength{\tabcolsep}{8pt}
    \centering
    \begin{tabular}{cccccccc}
    \toprule
    $C_2$ & $C^{\text{u}}_2$ & $C^{\text{d}}_2$ & F-SC cost & P-SC cost & $\Delta$cost & $q_0$ & $q_1$ \\
    \midrule
    15 & 20 & 10  & 418.59\euro & 416.91\euro & 0.40\% & -0.277 & 0.982 \\
    15 & 15 & 10  & 404.40\euro & 391.88\euro & 3.10\% & -0.253 & 0.899 \\
    15 & 20 & 15  & 413.65\euro & 412.93\euro & 0.17\% & -0.285 & 1.009 \\
    \bottomrule
    \end{tabular}
    \captionsetup{justification=centering}
    \caption{Cost savings in percentage under different values of $G_2$'s  balancing energy costs $C^{\text{u}}_2$ and $C^{\text{d}}_2$ (both given in \euro/MWh).}
    \label{table:toy_exp_1}
\end{table}

\begin{table}
    \setlength{\tabcolsep}{8pt}
    \centering
    \begin{tabular}{cccccccc}
    \toprule
    $C_2$ & $C^{\text{u}}_2$ & $C^{\text{d}}_2$ & F-SC cost & L-SC cost & $\Delta$cost & $q_0$ & $q_1$ \\
    \midrule
    15 & 20 & 10  & 418.59\euro & 444.05\euro & -6.08\%  &  4.967  & 0.964 \\
    15 & 15 & 10  & 404.40\euro & 418.84\euro & -3.57\%  &  5.490  & 0.883 \\
    15 & 20 & 15  & 413.65\euro & 680.92\euro & -64.61\% & 36.807  & 0.722 \\
    \bottomrule
    \end{tabular}
    \captionsetup{justification=centering}
    \caption{Cost savings in percentage under different values of $G_2$'s  balancing energy costs $C^{\text{u}}_2$ and $C^{\text{d}}_2$ (both given in \euro/MWh). Constraints \eqref{eq:estimation_problem_MO1}--\eqref{eq:estimation_problem_MOfinal} have been dropped from the training model~\eqref{eq:estimation_problem}.}
    \label{table:toy_exp_1_no_merit}
\end{table}

Therefore, due to the catastrophic impact that removing the merit-order constraints \eqref{eq:estimation_problem_MO1}--\eqref{eq:estimation_problem_MOfinal} from the training model~\eqref{eq:estimation_problem} may have on the actual performance of the obtained affine model, the strategy L-SC is no longer considered in the rest of our analysis.

\subsection{Impact of grid congestion}\label{sec:example:grid}
Here we introduce a variant of the base case in which the capacity of Line~1 has been set to 30 MW. Recall that the capacity of this line in the base case is unlimited, which we denote by symbol ``$\infty$" in Table~\ref{table:toy_exp_2}. The results collated in this new table are analogous to those in Table~\ref{table:toy_exp_1}. 

Recall that the estimation problem~\eqref{eq:estimation_problem}, whereby we determine the affine function $\widehat{L} = q_0 + q_1  L^{\rm F}$ , explicitly accounts for network constraints. In contrast, the computation of the net-demand point forecast $L^{\rm F}$ is typically based on statistical criteria alone and, consequently, ignores any possible limiting effect of the grid.

When the capacity of Line~1 is limited to 30 MW, our approach strongly corrects the point prediction $L^{\rm F}$ downwards, so that $\widehat{L}$ is kept in between 16 and 32 MW approximately. Thus, unit $G_1$ is dispatched well below the expected demand. This is clever because, in doing so, no (expensive) downward regulation from this unit has to be deployed in real time to comply with the limiting capacity of Line~1. In this way, the eventual realized demand at bus~3 can be satisfied, instead, with cheaper up-regulation from unit $G_2$ through Line~2. The ultimate result is that using $\widehat{L}$, given by our approach, in the forward problem~\eqref{eq:merit_order} is way more profitable than using the raw point forecast $L^{\rm F}$. 
\begin{table}
    \setlength{\tabcolsep}{8pt}
    \centering
    \begin{tabular}{ cccccc }
    \toprule
     $\overline{F}_1$ (MW)  & F-SC cost & P-SC cost & $\Delta$cost & $q_0$ & $q_1$ \\
    \midrule
    $\infty$ &  418.59\euro & 416.91\euro  &  0.40\% & -0.277 & 0.982 \\
    30       & 1034.70\euro & 724.46\euro  & 29.98\% & 15.725 & 0.175 \\
    \bottomrule
    \end{tabular}
    \captionsetup{justification=centering}
    \caption{Impact of grid congestion on cost savings.}
    \label{table:toy_exp_2}
\end{table}

\textcolor{black}{Based on this analysis, if network constraints were accounted for at the day-ahead stage, that is, in the optimization problem~\eqref{eq:merit_order}, the cost savings achieved by our approach would decrease, at the expense of significantly complicating the resolution of the bilevel problem~\eqref{eq:generic_binivel}.}

 

%
\subsection{Impact of the peak demand}
Now we change the peak demand and consider two variants of the base case in which we take $\overline{L} = 50$ MW and $\overline{L} = $ 150 MW (in the base case, $\overline{L} =$ 100 MW). The results of this new analysis are compiled in Table~\ref{table:toy_exp_3}. 

Again, as in the analysis of the impact of $G_2$'s balancing energy costs in Section~\ref{IE:reg_cost}, our approach systematically corrects the net-demand point forecast downwards to reduce the usage of down-regulation from $G_1$ in favor of the up-regulation from $G_2$. However, the cost savings achieved by our approach get diluted as the peak demand is augmented. The reason for this is twofold. First, the probability of events where the net demand takes on a value below the capacity of unit $G_1$ diminishes with growing $\overline{L}$. For instance, when $\overline{L}$ = 50 MW, the probability that the net demand is smaller than the capacity of $G_1$ is equal to one, which explains why our method delivers the highest cost savings in this variant (from among the three cases considered in this analysis). In contrast, as  $\overline{L}$ grows, that probability diminishes and the cheaper downward balancing generation from $G_2$ becomes more available. Second, the balancing costs account for a lower percentage of the total costs as the peak demand $\overline{L}$ increases.

%
%
\begin{table}
    \setlength{\tabcolsep}{12pt}
    \centering
    \begin{tabular}{ cccccc }
    \toprule
    $\overline{L}$ &  F-SC cost & P-SC cost & $\Delta$cost & $q_0$ & $q_1$ \\
    \midrule
    50  & 182.92\euro & 181.55\euro & 0.75\% & -0.138 & 0.982 \\
    100 & 418.59\euro & 416.91\euro & 0.40\% & -0.277 & 0.982 \\
    150 & 751.73\euro & 750.54\euro & 0.16\% & -0.421 & 0.997 \\
    \bottomrule
    \end{tabular}
    \captionsetup{justification=centering}
    \caption{Impact of peak demand.}
    \label{table:toy_exp_3}
\end{table}

 

%
\subsection{Impact of the net demand regime}
We conclude this small example by studying how the net demand regime affects the prescriptive power of the affine function $\widehat{L} = q_0 + q_1  L^{\rm F}$ that we determine by way of problem~\eqref{eq:estimation_problem}. To this end, we modify the support of the uniform distribution from which the per-unit net-demand point prediction  is randomly drawn. Thus, we distinguish a \emph{low-demand} regime, with $L^{\rm F} \sim \overline{L} \cdot {\rm U}(0.03,0.5)$, and a \emph{high-demand} regime, with $L^{\rm F} \sim \overline{L} \cdot {\rm U}(0.5,0.97)$. We also consider the base case, where $L^{\rm F} \sim \overline{L} \cdot {\rm U}(0.03,0.97)$ and therefore, no demand regime is differentiated. The corresponding results are provided in Table~\ref{table:toy_exp_4}.

In line with the observations in the previous analysis of the impact of the peak demand, under a low-demand regime, the expensive, but flexible unit $G_2$ is not dispatched in the forward problem. The downward correction to the net-demand point forecast our approach prescribes is then intended to benefit from the up-regulation provided by $G_2$, which is clearly more competitive than the down-regulation offered by $G_1$. The system features, therefore, a distinct cost asymmetry given by the expensive downward balancing generation of $G_1$ versus the cheap upward balancing generation of $G_2$. Our approach sees this asymmetry and corrects the net-demand point forecast downwards accordingly. In addition, since the beta distribution modeling the point forecast error is right-skewed for low levels of demand, said correction leads to substantial cost savings. In contrast, under a high-demand regime, $G_2$ is very likely to participate in the forward dispatch, whereas there is a lower probability that $G_1$ be needed to provide downward balancing energy, since the distribution of the point forecast error is left-skewed. Consequently, the cost structure of the system looks very different under a high-demand regime, which prompts a quite different affine function and reduces the cost savings obtained from our method.

Most importantly, in the base case, when no net-demand regime is distinguished, most of the benefits our approach can potentially bring for low values of net demand are lost. This motivates us to cluster net-demand observations into different regimes and use optimization problem~\eqref{eq:estimation_problem} to compute a possibly different affine model in the form $\widehat{L} = q_0 + q_1  L^{\rm F}$ for each demand regime, similarly to segmented regression in classical statistics. This is formalized in the next section.   
\begin{table}
    \setlength{\tabcolsep}{8pt}
    \centering
    \begin{tabular}{ cccccc }
    \toprule
    U(a, b) & F-SC cost & P-SC cost & $\Delta$cost & $q_0$ & $q_1$ \\
    \midrule
    U(0.03, 0.97) & 418.59\euro & 416.91\euro & 0.40\% & -0.277 & 0.982 \\
    U(0.03, 0.50) & 239.60\euro & 234.54\euro & 2.11\% & -0.102 & 0.917 \\
    U(0.50, 0.97) & 587.82\euro & 586.42\euro & 0.24\% & -6.646 & 1.088 \\
    \bottomrule
    \end{tabular}
    \captionsetup{justification=centering}
    \caption{Impact of the net-demand regime.}
    \label{table:toy_exp_4}
\end{table}

\section{Data clustering and partitioning}\label{sec:Part}
Take $\mathcal{N}:= \{1, \ldots, i, \ldots, N\}$, that is, the index set of the data sample $\{(\mathbf{x}_1, L_1)$, $\ldots$, $(\mathbf{x}_i, L_i)$, $\ldots$, $(\mathbf{x}_N, L_N)\}$ with $\mathbf{x}_i \in \mathbb{R}^p$ and $L_i \in \mathbb{R}^+, \forall i \in \mathcal{N}$. We partition $\mathcal{N}$ into a collection $\{\mathcal{N}_k\}_{k=1}^K$ of $K$ subsets that are pairwise disjoint and whose union is equal to $\mathcal{N}$. Consider the one-to-one mapping $\phi: \mathcal{N} \rightarrow \{1, 2, \ldots, K\}$, such that $\phi(i) = k$ if data point $(\mathbf{x}_i, L_i) \in \mathcal{N}_k$. Therefore, $\mathcal{N}_k:=\{i \in \mathcal{N}: \phi(i)=k\}$.

We compute $K$ affine models of the form  $\widehat{L} = \mathbf{q}^{\top}_{k}\mathbf{x}$, $k \leq K$, by solving the estimation problem~\eqref{eq:estimation_problem} for each subset sample $\mathcal{N}_k$. In practice, this means replacing $N$ and $\mathcal{N}$ in~\eqref{eq:estimation_problem} with $|\mathcal{N}_k|$  and $\mathcal{N}_k$, respectively.

To construct a meaningful mapping $\phi$, we employ the $K$-means algorithm that is implemented in the Python package \emph{scikit-learn} \citep{scikit-learn}, using the Euclidean distance. We note that, to construct $\phi$, this algorithm receives the feature sample $\{\mathbf{x}\}_{i \in \mathcal{N}}$ as input. \textcolor{black}{Once the mapping $\phi$ is computed, it is used to determine the cluster to which a new feature vector $\mathbf{x}$ belongs. That cluster is the one that minimizes the distance between its centroid and the feature vector $\mathbf{x}$.} That is, given a new observation of $\mathbf{x}$, say $\mathbf{x}_{N+1}$, $\phi(\mathbf{x}_{N+1}) = k$ means that $\mathbf{x}_{N+1}$ is predicted to belong to partition $\mathcal{N}_k$, and therefore, $\widehat{L} = \mathbf{q}^{\top}_k\mathbf{x}_{N+1}$ is to be used in the forward problem~\eqref{eq:merit_order}.

On a different issue, the estimation problem~\eqref{eq:estimation_problem} is a MIP program and, as such, computationally expensive in general. Actually, the size of~\eqref{eq:estimation_problem} grows linearly with the sample size. To keep the time to solve~\eqref{eq:estimation_problem} reasonably low, we reduce the cardinality of subsets $\{\mathcal{N}_k\}_{k=1}^K$ by means of the PAM K-medoids algorithm \citep{kaufman2009k-med} through the Python package implementation \emph{scikit-learn-extra}. This algorithm selects the most representative data points within each subset $\mathcal{N}_k$, the so-called \emph{medoids}, by minimizing the sum of distances between each point in $\mathcal{N}_k$ and said medoids. We remark that this reduction process results in data points (the medoids) with \emph{unequal} probability masses, so extra care should be taken when formulating objective function~\eqref{eq:estimation_problem_of} for each subset $\mathcal{N}_k$ considering the medoids only.  More specifically, the uniform weight $\frac{1}{N}$ appearing in the objective function~\eqref{eq:estimation_problem_of} should be replaced with a medoid-dependent weight representing the probability mass assigned to each medoid as a result of the reduction process.



\section{Case Study}\label{sec:CaseStudy}
In this section we assess the performance of our approach in a realistic case study that is based on the stylized model for the European power system that is described in~\cite{nahmmacher2014report}. Accordingly, we consider a pipeline network model with 28 nodes, each representing an European country. The capacities of the lines are also obtained from \cite{nahmmacher2014report}, in particular, we take the values from ``Table 14. Transmission capacities between model regions (GW)" that correspond to the year 2020. We assume that each node in the network (i.e., each European country) includes two types of power plants technologies, which we denote as \emph{base} and \emph{peak}, respectively. Again, the available capacity of both technologies has been assigned based on the data in~\cite{nahmmacher2014report} corresponding to 2020 for each country. More specifically, the base power-plant capacities have been obtained by adding up the installed capacities of the technologies ``Nuclear'', ``Hard coal'', ``Oil'' and ``Lignite'' and the peak power-plant capacities from the technologies ``Natural Gas'', ``Waste'' and ``Other gases''. The nodes of the system and the resulting generation capacities of each type are listed in Table~\ref{tab:base_peak_node_cap}.

\begin{table}
\setlength{\tabcolsep}{4pt}
{\footnotesize
    \begin{tabular}{lllllllllllllll}
    \hline
    Country & AT  & BE  & BG   & CH  & CZ   & DE   & DK   & EE  & ES   & FI  & FR   & GB   & GR  & HR  \\
    \hline
    base    & 0.4 & 6.1 & 6.7  & 3.4 & 14.4 & 46   & 2.4  & 2   & 16.3 & 6.5 & 68.3 & 20.6 & 3.9 & 1.3 \\
    peak    & 5.9 & 6.8 & 1    & 0.6 & 1.3  & 27.9 & 1.7  & 0.2 & 29.6 & 3.6 & 11.9 & 31.2 & 4.9 & 0.7 \\
    \hline
    Country & HU  & IE  & IT   & LT  & LU   & LV   & NL   & NO  & PL   & PT  & RO   & SE   & SI  & SK  \\
    \hline
    base    & 3.3 & 1.8 & 8.7  & 0   & 0    & 0    & 4.5  & 0   & 27.8 & 1.8 & 5.4  & 11.1 & 1.8 & 2.7 \\
    peak    & 4.1 & 4.2 & 46.2 & 1.8 & 0.1  & 1.2  & 19.3 & 0   & 3.5  & 4.6 & 2.9  & 1.1  & 0.7 & 1.5\\
    \hline
    \end{tabular}
}
    \caption{Base and peak generation capacity (GW) installed per node of the European network.}
    \label{tab:base_peak_node_cap}
\end{table}

To build a data sample of the form $\{(\mathbf{x}_i,L_{i})\}_{i \in \mathcal{N}}$, we have collected the actual aggregate hourly demand, wind, solar and hydro energy production for each country (node of the system) in 2020 from the ENTSO-e Transparency Platform~\citep{Entsoe}. We have also retrieved the day-ahead forecast of the hourly demand and the produced wind and solar energy from this platform. To get series of net demand values (both forecast and actual), we have subtracted the respective wind, solar and hydro power data series from the aggregate day-ahead forecast/actual demand series. We clarify that no day-ahead forecast for the hydro power production is available in~\cite{Entsoe}, so the series of real hydro power production has been used (instead of the missing day-ahead hydro forecast) for the computation of day-ahead forecasts of the nodal net demands.  Some minor gaps in the data extracted from \cite{Entsoe} have been filled through linear interpolation. 

The marginal costs of energy generation and up- and downward balancing energy provision of each unit are randomly sampled from the uniform distributions specified in Table~\ref{tab:case_study_prices}. The so-obtained values for these costs have remained fixed throughout the experiments performed in this section. We point out that in the uniform distributions of Table~\ref{tab:case_study_prices}, we have considered that base power units are cheap but inflexible, and thus, with costly balancing energy. In contrast, peak power plants are expensive, but flexible, and hence, with more competitive balancing energy costs.

\begin{table}
    \centering
    \begin{tabular}{cccc}
        \hline
              &   $C$    &   $C^u$   &   $C^d$   \\
        \hline      
        base  & U(8, 12)  & U(60, 70) & U(-40, -50) \\
        peak  & U(36, 44) & U(45, 50) & U(30, 35)   \\
        \hline
    \end{tabular}
    \caption{Uniform distributions from which the marginal production, up- and down-balancing costs (in \euro/MWh) of the units in the European system have been sampled.}
    \label{tab:case_study_prices}
\end{table}

We conduct a rolling simulation on the data of 2020, in which we gradually select non-overlapping windows of 150 points each. From each window, we randomly sub-sample (without replacement) the indexes corresponding to the training and test sets, which are eventually made up of 100 and 50 samples, respectively. We take ten windows over which we average the results that follow. \textcolor{black}{First, we compute the average cost of the baseline method F-SC, i.e., the method that uses the point forecast $L^{\rm F}$ of the net demand as $\widehat{L}$ in~\eqref{eq:merit_order}, which amounts to 5092.8k\euro. For comparison, we also compute the average cost of a perfect information benchmark in which $\widehat{L}$ is equal to the actual realization of the net demand $L$. This perfect forecast approach is unrealistic but can be used to assess how much could be gained from improving the prescribed net demand. For this case study, the perfect information cost is 3380.9k\euro, i.e, 33.6\% lower than that delivered by method F-SC.}

As in the example of Section~\ref{sec:Ex}, we consider a feature vector $\mathbf{x}$ made up of the day-ahead forecast of the system net demand, $L^{\rm F}$, (measured in MWh), enlarged with an additional feature fixed to one to accommodate the intercept of the affine models $\widehat{L} = \mathbf{q}^{\top}_k\mathbf{x} = q_{0k} + q_{1k} L^{\rm F} $, $k \leq K$. 


In the analysis we conduct next, we consider various values for $K$ (number of partitions and hence of affine models) and several percentage reductions of the number of data in each partition $\mathcal{N}_k, k \leq K$. The results of this analysis are summarized in Table~\ref{tab:med_seg_results}, where ``$r$\%'' in the first column 
means that only that percentage of medoids in the partition $\mathcal{N}_k$ (more precisely $\ceil*{\frac{r}{100} |\mathcal{N}_k|}$, where $\lceil \cdot \rceil$ denotes the ceiling operator) have been used to estimate the affine function $\widehat{L} = \mathbf{q}^{\top}_k\mathbf{x}$ through~\eqref{eq:estimation_problem}. 
\textcolor{black}{This table shows the average cost achieved by our approach for different values of $r$ and $K$, the cost savings with respect to the cost of the baseline method F-SC (5092.8k\euro), and the average time the solution to the K estimation problems~\eqref{eq:estimation_problem} takes.} The reported cost savings have been computed out of sample, that is, on the test sets. Beyond the fact that these savings are significant in general, it is clear that our prescriptive approach benefits from exploiting different affine models under different net-demand regimes, which confirms the preliminary conclusion we draw in this regard through the small example of Section~\ref{sec:Ex}. Nevertheless, it is also true that the added benefit rapidly plateaus as $K$ grows. Actually, the bulk of the economic gains we get through the partitioning of the data sample is already reaped with $K = 2$. On the other hand, increasing $K$ has a positive side effect: It remarkably reduces the time to solve the MIP problem~\eqref{eq:estimation_problem}. In addition, this time can be shortened even further, with a tolerable reduction in cost savings, by using only the medoids of the partitions $\mathcal{N}_k, k \leq K$,  when estimating the affine models through~\eqref{eq:estimation_problem}. {\color{black} Notwithstanding this, $K$ cannot be made arbitrarily big since it can lead to \emph{overfitting}. Therefore, the choice of $K $ should be guided by a training-validation scheme similar to the one we show in Table~\ref{tab:med_seg_results} to ensure that increasing $K$ does not lead to a reduction in the benefits of the affine rules out of sample.}

\begin{table}[]
    \centering
    \begin{tabular}{ccccc}
    \hline
    $r$ & $K$ & P-SC cost & $\Delta$cost & Time (s) \\  
    \hline
    100\% & 1 & 4948.7k\euro & 2.83\% & 2127.7 \\
    100\% & 2 & 4874.4k\euro & 4.29\% & 283.7 \\
    100\% & 5 & 4851.6k\euro & 4.74\% & 75.9 \\
    100\% & 7 & 4851.1k\euro & 4.75\% & 28.0 \\
    \hline
    50\% & 1 & 4956.9k\euro & 2.67\% & 180.0 \\
    50\% & 2 & 4877.3k\euro & 4.23\% & 27.2 \\
    50\% & 5 & 4869.4k\euro & 4.39\% & 7.4 \\
    50\% & 7 & 4886.1k\euro & 4.06\% & 5.5 \\
    \hline
    20\% & 1 & 4971.6k\euro & 2.38\% & 8.3 \\
    20\% & 2 & 4883.2k\euro & 4.12\% & 3.2 \\
    20\% & 5 & 4882.8k\euro & 4.12\% & 1.1 \\
    20\% & 7 & 4890.9k\euro & 3.97\% & 1.4 \\
    \hline
    \end{tabular}
    \caption{Average cost savings and average time to solve the estimation problem~\eqref{eq:estimation_problem} for a number $K$ of partitions and various levels $r$ of reduction in the size of the original training sets (in percentage).}
    \label{tab:med_seg_results}
\end{table}

\begin{figure} 
\centering
\begin{tikzpicture}[scale=0.9, every node/.style={scale=0.9}]
	\begin{axis}[
	width=0.7\textwidth,
    xmin=150, xmax=300, ymin=150, ymax=300,
    xtick = {150, 200, 250, 300},
    ytick = {150, 200, 250, 300},
	legend style={at={(0.5,0.85)},anchor=south,legend cell align=center,legend columns=5},	
	clip marker paths=true,	
	xlabel = Predicted aggregate net demand ($L^{\rm F}$),
	ylabel = Prescribed aggregate net demand $\left(\widehat{L}\right)$,
	]
	\addplot[only marks, mark=square*, draw=black,fill=black] table [x=c900-1050_s5_m1.0_c1_FO, y=c900-1050_s5_m1.0_c1_CED, col sep=comma] {data/s5_m1_900-1050.csv}; \addlegendentry{$\mathcal{N}_1$}
	\addplot[only marks, mark=triangle*, draw=gray, fill=gray] table [x=c900-1050_s5_m1.0_c4_FO, y=c900-1050_s5_m1.0_c4_CED, col sep=comma] {data/s5_m1_900-1050.csv}; \addlegendentry{$\mathcal{N}_2$}
	\addplot[only marks,mark=star, draw=black] table [x=c900-1050_s5_m1.0_c5_FO, y=c900-1050_s5_m1.0_c5_CED, col sep=comma] {data/s5_m1_900-1050.csv}; \addlegendentry{$\mathcal{N}_3$}
	\addplot[only marks, mark=diamond*, draw=gray, fill=gray] table [x=c900-1050_s5_m1.0_c2_FO, y=c900-1050_s5_m1.0_c2_CED, col sep=comma] {data/s5_m1_900-1050.csv}; \addlegendentry{$\mathcal{N}_4$}
	\addplot[only marks, mark=*, draw=black,fill=black] table [x=c900-1050_s5_m1.0_c3_FO, y=c900-1050_s5_m1.0_c3_CED, col sep=comma] {data/s5_m1_900-1050.csv}; \addlegendentry{$\mathcal{N}_5$}
	\addplot[mark=none, dotted, domain=120:320, line width=0.75pt, draw=gray] {1 * x};
	\end{axis}	
\end{tikzpicture} \\
\caption{Prescribed affine transformation of the day-ahead net-demand forecast (aggregated system-wise). Demand is given in GW.}
\label{fig:case_study_991ada}
\end{figure}
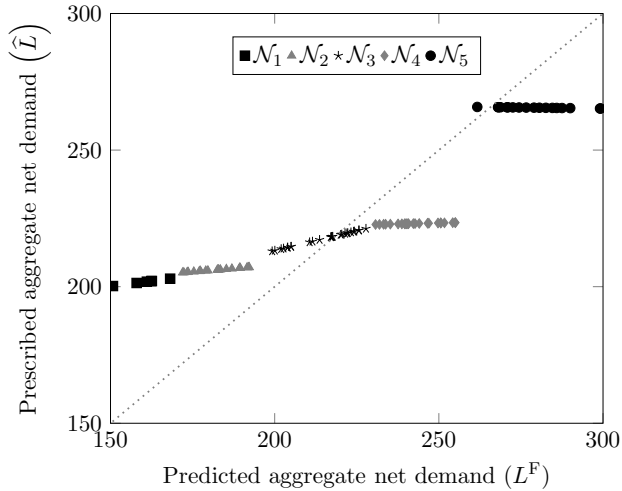

To comprehend where those cost savings our approach yields come from, in Figure~\ref{fig:case_study_991ada} we plot the predicted aggregate net demand $L^{\rm F}$ against the one prescribed by our method, i.e., $\widehat{L}$. The plot corresponds to one window of 150 data points taken at random out of the ten we have considered in the rolling-window simulation. Furthermore, the figure depicts results from the case with five partitions ($K = 5$). It can be seen that, when the system net demand is predicted low, our method prescribes to overestimate it. This prescription is motivated by two facts. On the one hand, the overestimation of the net demand in the forward problem is covered by cheap power plants, whereas it reduces the need for upward balancing. On the other, even though it slightly increases the demand for downward balancing, the group of units that down-regulate remains the same in any case, i.e., with and without the overestimation, due to the limitations of the network. As a result, the cost savings linked to the reduction in up-balancing outweigh the extra costs incurred by the increase in down-balancing. It is interesting to note that, as mentioned at the end of Section~\ref{sec:PD}, system operators, based on their accumulated experience, often introduce an upward bias into the net demand forecast \citep{CAISO}, precisely as our model here suggests.

As the level of net demand grows, the overestimation of the system net load that our method prescribes diminishes to a point where the prescribed amount flattens (see partitions $\mathcal{N}_4$ and $\mathcal{N}_5$). Again, this phenomenon is caused by the network and the limitations it imposes. Indeed, our method avoids dispatching power plants in the forward problem which, despite being their turn in the cost-merit order, would have been irretrievably down-regulated in real time because of network bottlenecks. For instance, in partition $\mathcal{N}_4$, F-SC consistently dispatches the DE base generator, with its massive 46 GW, to maximum capacity. However, due to grid constraints, this unit is subsequently down-redispatched to around 30 GW. On the contrary, P-SC takes into consideration that this power plant is one of the latest to be scheduled in this partition and foresees the grid limitation on the power flow, thus constraining the aggregated energy production and systematically dispatching such a unit to the previously mentioned 30 GW.

\textcolor{black}{Finally, Figure \ref{fig:case_study_991ada} also illustrates that, even though the prescriptive functions $h_{\boldsymbol{q}}(\cdot)$ we employ are affine, the use of several clusters leads to a \emph{piece-wise} affine function that is able to approximate the non-linear relationships between the features and the prescribed parameters.}

\section{Conclusions}\label{sec:Conclu}

In this paper, we have proposed a data-driven method to prescribe the value of net demand that the forward stage in a two-stage generation scheduling problem should use in order to minimize the expected total cost of operating the underlying power system. For this purpose, we have formulated a mixed-integer linear program that trains an affine function to map the \emph{predicted} net demand into the \emph{prescribed} one. 

Numerical experiments conducted out of sample on a stylized model of the European electricity grid reveal that the cost savings implied by the estimated affine mappings are substantial, well above 2\%. Furthermore, on the grounds that the cost structure of a power system is highly dependent on its operating point, and hence, on the level of net demand, we have devised a $K$-means-based partition strategy of the data sample to train different affine mappings for different net-demand regimes. The utilization of this strategy is shown to have a positive twofold effect in the form of substantially increased costs savings and a remarkable drop in the computational burden of the proposed MIP training model. Finally, we have further complemented the partitioning of the data sample with a medoid-based reduction in the size of the partitions, achieving additional speedups in solution times. All this together opens up the possibility to leverage our prescriptive approach in larger instances.

Future work will include attempts to optimize the partitioning of the data sample by embedding it into the MIP training model {\color{black} and to devise a procedure to efficiently update the parameters of the affine rules as new data become available without having to solve a new instance of the estimation problem~\eqref{eq:generic_binivel}}. \textcolor{black}{Another direction for future research is the extension of the proposed methodology to prescribe net demand trajectories (of 24 hours, for example) for generation scheduling problems that include inter-temporal constraints such as ramping limits or minimum times.}



\section*{CRediT authorship contribution statement}

\textbf{J. M. Morales}: Conceptualization, Methodology, Investigation, Formal analysis, Writing - original draft, Supervision, Funding acquisition. \textbf{M. A. Muñoz}: Data curation, Software, Methodology, Investigation, Formal analysis, Validation, Writing - Original Draft. \textbf{S. Pineda}: Conceptualization, Methodology, Investigation, Formal analysis, Writing - original draft, Supervision.

\section*{Acknowledgments}

This work was supported in part by the European Research Council (ERC) under the EU Horizon 2020 research and innovation program (grant agreement No. 755705), in part by the Spanish Ministry of Science and Innovation (AEI/10.13039/501100011033) through project PID2020-115460GB-I00, and in part by the Junta de Andalucía (JA) and the European Regional Development Fund (FEDER) through the research project P20\_00153. M. Á. Mu\~{n}oz is also funded by the Spanish Ministry of Science and Innovation through the State Training Subprogram 2018 of the State Program for the Promotion of Talent and its Employability in R\&D\&I, within the framework of the State Plan for Scientific and Technical Research and Innovation 2017-2020 and by the European Social Fund (reference PRE2018-083722 and ENE2017-83775-P). Finally, the authors thankfully acknowledge the computer resources, technical expertise, and assistance provided by the SCBI (Supercomputing and Bioinformatics) center of the University of M\'alaga.


\bibliographystyle{elsarticle-harv}
\bibliography{bibliography}

\end{document}